\documentclass[a4paper,12pt]{article}
\usepackage[dvips]{epsfig}
\usepackage{amsmath,amssymb,euscript,graphicx,amsfonts}
\usepackage{enumerate,color}
\usepackage{graphicx}
\usepackage{hyperref}
\usepackage{tikz}

\newtheorem{theorem}{Theorem}[section]
\newtheorem{proposition}[theorem]{Proposition}
\newtheorem{lemma}[theorem]{Lemma}
\newtheorem{corollary}[theorem]{Corollary}

\newtheorem{problem}[theorem]{Problem}

\newtheorem{remark}[theorem]{Remark}

\newcommand{\Qed}{\rule{2.5mm}{3mm}}

\renewcommand{\P}{\mathcal{P}}
\newcommand{\Q}{\mathcal{Q}}

\newcommand{\B}{\mathcal{B}}

\newcommand{\ZZ}{\mathbb{Z}}

\newcommand{\FF}{\mathbb{F}}

\newcounter{case}

\renewcommand{\thecase}{\arabic{case}}

\newcounter{subcase}

\numberwithin{subcase}{case}

\def\syl{\hbox{\rm Syl}}
\def\soc{\hbox{\rm Soc}}

\def\PGL{\hbox{\rm PGL}}

\def\SL{\hbox{\rm SL}}

\def\PSL{\hbox{\rm PSL}}
\def\AGL{\hbox{\rm AGL}}
\def\char{\hbox{\rm char}}
\def\demo{{\bf Proof}\hskip10pt}
\def\di{\bigm|} \def\lg{\langle} \def\rg{\rangle}
\def\qqed{\hfill $\Box$}
\def\a{\alpha}    
\def\r{\rho}

\def\Aut{\hbox{\rm Aut\,}}

\def\a{\alpha}    
\def\r{\rho}

 \def\og{\overline G} \def\oh{\overline H}

   \def\o1{\overline 1}

\def\o{\overline}   
\def\di{\bigm|} \def\lg{\langle} \def\rg{\rangle}

\def\di{\bigm|} \def\lg{\langle} \def\rg{\rangle}

\usepackage{minitoc}

\begin{document}

	\begin{center}
		{\bf\large  On Hamilton cycles in connected vertex-transitive graphs of order 
			$2pq$
			}
		\footnote{This work was supported in part  by the National Natural Science Foundation of China
			(12471332).}
	\end{center}


	\begin{center}
	 Tianlei Zhou\footnote{Corresponding author: 2230501007@cnu.edu.cn}\\
		\medskip
		{\it {\small
				School of Mathematical Sciences, Capital Normal University,\\
				Beijing 100048, People's Republic of China\\
		}}
	\end{center}

	\renewcommand{\thefootnote}{\empty}

\begin{abstract}
The existence of Hamilton cycles in connected vertex-transitive graphs is an open problem in algebraic graph theory, closely related to Lov\'asz's 1969 conjecture on Hamilton paths. All connected vertex-transitive graphs of order $pq$, where $p$ and $q$ are distinct primes, are known to be Hamiltonian except for the Petersen graph, and the primitive case for graphs of order $2pq$ has been resolved except for the Coxeter graph. This paper considers connected non-quasiprimitive graphs of order $2pq$ admitting a transitive subgroup of the automorphism group with an intransitive normal subgroup inducing prime-length orbits. We prove that all such graphs contain a Hamilton cycle. This result covers a class of non-quasiprimitive graphs of order $2pq$, advancing the full resolution of the $2pq$ case.
\end{abstract}

\section{Introduction}
\label{sec:intro}
\indent

The existence of Hamilton cycles and Hamilton paths in vertex-transitive graphs is one of the most long-standing open problems in algebraic graph theory, originating from Lov\'asz's 1969 conjecture\cite{L70} which asks whether every finite connected vertex-transitive graph contains a Hamilton path. In 1981, Alspach  \cite{A83} asked if there exists an infinite number of connected vertex-transitive graphs that do not have a Hamilton cycle. Over half a century of intensive study, no connected vertex-transitive graph lacking a Hamilton path has been found, and only four connected vertex-transitive graphs on at least three vertices are known to not contain a Hamilton cycle: the Petersen graph, the Coxeter graph, and the two graphs obtained by replacing each vertex of these two with a triangle. Notably, none of these four exceptional graphs is a Cayley graph, leading to the folklore conjecture that every connected Cayley graph on at least three vertices contains a Hamilton cycle, see  \cite{CG96,D83, GWM14,GM07,GKMM12,KM09,MR234, DM83,WM18,WM15}. While the full conjecture remains open, substantial progress has been made by resolving the Hamiltonian problem for vertex-transitive graphs of specific orders, often relying on the classification of finite simple groups and detailed analysis of permutation group actions. Chen investigated a stronger structural property in 1988 (see \cite{chen}). A Hamiltonian graph is called edge-Hamiltonian if every edge lies on some Hamilton cycle. Chen proposed the following conjecture which asks whether every Hamiltonian Cayley graph is edge-Hamiltonian.

This conjecture has been verified for several classes of Cayley graphs, including those on abelian groups, groups of order \(2^k\), and dihedral groups\cite{chen}. Since Cayley graphs constitute a natural subclass of vertex-transitive graphs, one may ask whether the same phenomenon extends to all vertex-transitive graphs. This motivates the following conjecture.
\begin{problem}
If a vertex-transitive graph is Hamiltonian, then every edge lies in some Hamilton cycle.
\end{problem}
\vskip 3mm

Early work in this direction established the existence of Hamilton cycles (or paths) for vertex-transitive graphs of orders $kp$ for \(k\leq 10\), \(p^j\) for \(j\leq 5\), and \(2p^2\) for prime $p$, see \cite{C98,MR4328721,MR4548744,KM08,KMZ12,KS09,DM87,DM92,DM85,MP82,MP83,MS4,MS2,MS5,Z15} and a survey paper \cite{KM09}. A landmark breakthrough came in 2021, when Du, Kutnar and Marušič completely resolved the Hamiltonian problem for vertex-transitive graphs of order $pq$, where $p$ and $q$ are distinct primes: they proved that with the sole exception of the Petersen graph, every connected vertex-transitive graph of order $pq$ contains a Hamilton cycle (see \cite{MR4328721}). This result closed the problem for two-prime orders, and naturally directed attention to the next natural case: vertex-transitive graphs of order $2pq$, where $p$ and $q$ are distinct primes.

\vskip 3mm
The $2pq$ order case is significantly more complex, as its full automorphism group has more types and the group action is more complex. Du, Tian and Yu showed that  except for the Coxeter graph, every graph of such order  contains a Hamilton cycle, provided the automorphism group acts primitively on its vertices (see \cite{MR4548744}). Most recently, Du and Zhou proved that every connected vertex-transitive graph of order $6p$ (corresponding to \(q=3\)) contains a Hamilton cycle, with the only exception being the triangle-replaced Petersen graph.
\vskip 3mm

The main result of this paper is the following theorem.

\begin{theorem}
	\label{the:main}
Every connected non-quasiprimitive graph of order $2pq$ admitting a transitive subgroup of its automorphism group with an intransitive normal subgroup inducing prime-length orbits contains a Hamilton cycle, where $p$ and $q$ are distinct primes.
\end{theorem}

This paper is organized as follows: after this introductory section, some preliminaries will be given in Section 2 and the proof of Theorem~\ref{the:main} will be given in Section 3.
\section{Terminology, notation and some useful results}
\label{sec:pre}
This section consists of four subsections dealing with basic definitions and notation, coset graphs, existence of Hamilton cycles in particular graphs, and group- and   finite field-theoretic results, respectively.

\subsection{Basic definitions and notation}
\label{ssec:definition}
\noindent Throughout this paper, graphs are finite, simple and undirected, and groups are finite. Given a graph $\Gamma$, we denote the vertex set by $V(\Gamma)$ and the edge set by $E(\Gamma)$, respectively. For two adjacent vertices $u,v \in V(\Gamma)$, we write $u \sim v$ and denote the corresponding edge by $\{u,v\}$. For $\a\in V(\Gamma)$, we denote the neighborhood of $\a$ by $N(\a)$. Let $U$ and $W$ be two disjoint subsets of $V(\Gamma)$. Then by $X( U )$ and $X[U,W]$ we denote the subgraph of $\Gamma$ induced by $U$ and the bipartite subgraph of $\Gamma$ induced by the edges having one end-vertex in $U$ and the other end-vertex in $W$, respectively.

\vskip 3mm
Let $G$ be a group acting faithfully and transitively on a set $V$. A nonempty subset $B$ of $V$ is called a $block$ for $G$ if for each $g\in G$ either $B^g=B$ or $B^g\cap B=\emptyset$, where $V$ and the singletons $\{v\}$ $(v\in V)$ are called the {\it trivial} blocks, and other blocks are called {\it nontrivial} blocks. Put $\B=\{B^g: g\in G\}$. Then the sets in $\B$ form a partition of $V$. We call $\B$ the {\em system of blocks} containing $B$. If $G$ has no nontrivial blocks on $V$, then $G$ is {\em primitive}, and is {\em imprimitive} otherwise. If every nontrivial normal subgroup of $G$ is transitive on $V$, then $G$ is {\em quasiprimitive}, and is {\em non-quasiprimitive} otherwise. If $G$ acts non-quasiprimitively on $V$ with an intransitive normal subgroup $N$, then for any   $v\in V$, $v^N=\{v^g: g\in N\}$ is called an $N$-{\em block}. $|v^N|$ is called the length of blocks. The set of $N$-blocks forms a system of blocks for $G$.
\vskip 3mm
A graph $\Gamma$ is said to be {\it vertex-transitive} if its full automorphism group  $\Aut(\Gamma)$ acts transitively on $V(\Gamma)$. A vertex-transitive graph $\Gamma$ is said to be {\it quasiprimitive} if there exists a group $G\le \Aut(\Gamma)$ such that $G$ acts quasiprimitively on $V(\Gamma)$ and is {\it non-quasiprimitive} otherwise. In particular, if there exists a group $G\le\Aut(\Gamma)$ such that $G$ acts primitively on $V(\Gamma)$, then $\Gamma$ is said to be  {\it primitive} and {\it imprimitive} otherwise.
\vskip 3mm
Let $\Gamma$ be a non-quasiprimitive graph and $G$ be a transitive subgroup of $\Aut(\Gamma)$. Then $G$ contains an intransitive normal subgroup $N$. Let $\B$ be the set of $N$-blocks and $A,B\in \B$ be distinct. By $d(A)$ and $d(A,B)$ we denote the valency of $X(A)$ and $X[A,B]$, respectively. Since $N$ is transitive on both $V(X(A))$ and $V(X(B))$, we know that $X[A,B]$ is regular or is a null graph. A {\it normal quotient graph corresponding to} $N$ is the graph $\Gamma_\B$ whose vertex set is $\B$ with $A,B\in \B$ adjacent if there exist vertices $a \in A$ and $b \in B$, where $a,b\in V(\Gamma)$, such that $a \sim b$ in $\Gamma$. Let $N_0$ be the kernel of the action of $G$ on $\B$. $\Gamma_\B$ is also called the {\it $N$-block graph}. Then $N\le N_0$ and there exists a natural faithful action of $G/N_0$ on $\B$. In particular, the action is transitive. Thus, if $\Gamma$ is a connected non-quasiprimitive graph, then $\Gamma_\B$ is a connected vertex-transitive graph.
\vskip 3mm
Let $m\geq 1$ and $n\geq 2$ be integers. An automorphism $\rho$ of a graph $\Gamma$ of order $mn$ is called $(m,n)$-{\em semiregular} (in short, {\em semiregular}) if as a permutation on $V(\Gamma)$ it has a cycle decomposition consisting of $m$ cycles of length $n$. The question whether all vertex-transitive graphs admit a semiregular automorphism  is one of the famous open problems in algebraic graph theory (see, for example, \cite{bcc15,seven,DMMN07,G1,M81}).
\vskip 3mm
Let $\P$ be the set of orbits of $\lg \rho\rg$. Let $A, B \in \P$ be distinct. Similarly, the graph $X[A,B]$ is regular or is a null graph. Let the {\em quotient graph corresponding to $\P$} be the graph $X_\P$ whose vertex set is  $\P$ with $A, B \in \P$ adjacent if $X[A,B]$ is not a null graph.

\vskip 3mm

A graph $\Gamma$ is called an $(m,n)$-{\em metacirculant}, where $m$ and $n$ are positive integers, if $\Gamma$ is of order $mn$ and has two automorphisms $\rho$, $\sigma$ such that

\begin{enumerate}
\item  $\rho$ is $(m,n)$-semiregular;

\item  $\sigma$ cyclically permutes the $m$ orbits of $\lg \rho\rg$ and normalizes $\lg \rho\rg$; and

\item  $\sigma^m$ fixes at least one vertex of $\Gamma$.
\end{enumerate}

It follows that $\lg\rho,\sigma\rg$ is a vertex-transitive subgroup of $\Aut(\Gamma)$ isomorphic to $\ZZ_n\rtimes\ZZ_t$, where $t=o(\sigma)$ and $m\mid t$.

\subsection{Coset graphs and minimal connected coset graphs}
\label{ssec:orbital}
 \subsubsection{Coset graphs}
Let $\Gamma$ be a connected vertex-transitive graph, and let $G$ be a transitive subgroup of $\Aut(\Gamma)$. Fix an arbitrary vertex \(\alpha \in V(\Gamma)\). Write \(H = G_\alpha = \{g \in G : \alpha^g = \alpha\}\) for the stabilizer of \(\alpha\) in $G$. There is a natural $G$-equivariant bijection between the vertex set \(V(\Gamma)\) and the right coset space \([G:H] = \{Hg : g \in G\}\), which maps each vertex \(\alpha^g\) to the coset $Hg$ for all \(g \in G\). Under this identification, the transitive action of $G$ on vertices corresponds exactly to the right multiplication action of $G$ on cosets. Let $\Omega=\{g:\a^g\in N(\a)\}$ and $D=\cup_{g\in\Omega}HgH$. Then $\Gamma$ is determined by $G,H$ and $D$. This correspondence allows a purely group-theoretic construction of vertex-transitive graphs. The {\it coset graph} (or {\it double coset graph}) \(X(G, H, D)\) is defined to have vertex set \([G:H]\), with adjacency given by: \(Hg_1 \sim Hg_2 \quad \text{if and only if} \quad Hg_2 g_1^{-1} \subseteq D.\) We denote the number of right cosets contained in $D$ by $|D|$, that is, $|D|=|\{Hg:g\in D\}|$. Similarly, $|Hg_0H|=|\{Hg:g\in Hg_0H\}|$ is called the length of the double coset $Hg_0H$.
\vskip 3mm
For $\Gamma=X(G,H,D)$ to be a finite, simple, undirected graph, three standard conditions hold:
\begin{enumerate}
\item \(H \cap D = \emptyset\) (no loops);
\item  \(D = D^{-1} = \{d^{-1} : d \in D\}\) (symmetric adjacency for undirected edges);
\item   $\Gamma$ is connected if $\lg g:g\in D\rg=G$.
\end{enumerate}
In particular, if $D_g=HgH=Hg^{-1}H$, then $D_g$ is said to be self-paired and is nonself-paired otherwise. If $D_g\neq D_{g^{-1}}$, then  $D_g$ and $D_{g^{-1}}$ are said to be symmetric-paired.

\vskip 3mm
$\Gamma=X(G,H,D)$ is $G$-symmetric provided $D=HgH=Hg^{-1}H$. $\Gamma=X(G,H,D)$ is $G$-edge transitive provided $D=HgH\cup Hg^{-1}H$, where $HgH\neq Hg^{-1}H$.

\vskip 3mm

Let $\Gamma=X(G,H,D)$ be a connected vertex-transitive graph and $N$ be a nontrivial intransitive normal subgroup of $G$. Suppose $|V(\Gamma)|=n$. Then $N$ induces blocks of length $m$, where $m\ge 2$ and $m \mid n$. Let $\B$ be the set of $N$-blocks and $A, B \in \cal{B}$ be two distinct blocks. Let $N_0$ be the kernel of the action of $G$ on $\B$. Then $\Gamma_\B\cong X(\og, \oh, \o{D})$, where $\overline{G}=G/N_0, \overline{H}=(HN_0)/N_0$ and $\overline{D}=\cup_{g\in D}\overline{HgH}\backslash {\overline{H}}$.

\subsubsection{Minimal connected coset graphs}

Let $\Gamma=X(G,H,D)$ be a coset graph, where $D=\cup_{g\in \Omega}HgH$.  Let $\Gamma_{0}=X(G,H,D_{0})$ be a subgraph of $\Gamma$ and $D_{g_{0}}=Hg_0H\cup Hg_0^{-1}H$, where $g_0\in D$ and $D_{0}=D\backslash D_{g_0}$. If $\Gamma$ is connected and, for any $g_0\in D$, $\Gamma_0$ is not connected, then $\Gamma$ is called a {\it $G$-minimal connected (coset) graph} (in short, {\it $G$-minimal graph}). That is, if $\Gamma$ is a $G$-minimal graph, then $\Gamma_{0}=X(G,H,D_0)$ is not connected, where $D_0\subsetneq D$ and $D_0=D_0^{-1}$.
\vskip 3mm
Every connected vertex-transitive graph $\Gamma=X(G,H,D)$ contains at least one $G$-minimal graph. If every  $G$-minimal graph contains a Hamilton cycle, then every connected vertex-transitive graph with a transitive group $G$ contains a Hamilton cycle. Thus, to find a Hamilton cycle in $\Gamma$, we only need to consider $G$-minimal graphs.

\subsection{Existence of Hamilton cycles in particular graphs}
\label{ssec:numbers}
The following  known results about the existence of Hamilton cycles in particular graphs will be used later.

\begin{proposition}
	\label{pro2.1}
{\rm \cite[Theorem~6]{BJ78}} Every $2$-connected regular graph of order $n$ and valency at least $n/3$ contains a Hamilton cycle.

\end{proposition}

\begin{proposition}
	\label{pro2.2}
{\rm \cite[Lemma~5]{MP82}} Let $X$ be a graph admitting an $(m, p)$-semiregular automorphism $\rho$, where $p$ is a prime. Let $C$ be a cycle of length $k$ in the quotient graph $X_{\mathcal{P}}$, where $\mathcal{P}$ is the set of orbits of $\rho$. Then the lift of $C$ either contains a cycle of length $kp$ or it consists of $p$ disjoint $k$-cycles. In the latter case, we have $d(S,S\sp{'})=1$ for every edge $SS\sp{ '}$ of $C$.
\end{proposition}

\begin{corollary}
	\label{co2.3}
Let $\Gamma=X(G,H,D)$ be a non-quasiprimitive graph. Let $N$ be an intransitive normal subgroup of $G$ and $\B$ be the set of $N$-blocks. We may assume that $N$ is the kernel of the action of $G$ on $\B$.  Suppose that $N$ induces blocks of length $r$, where $r$ is a prime. If the action of $N$ on each $N$-block is unfaithful and $\Gamma_\B$ contains a Hamilton cycle, then $\Gamma$ contains a Hamilton cycle.
\end{corollary}

By Corollary~\ref{co2.3}, to find a Hamilton cycle in $\Gamma$, if $\Gamma_\B$ contains a Hamilton cycle, we only need to consider the case where $N$ acts faithfully on $N$-blocks, which implies $N\lessapprox S_r$.

\begin{proposition}
	\label{pro2.4}
{\rm \cite[Theorem~1.2]{MR234}} Let $G$ be a finite group.  Every connected Cayley graph on $G$ with at least three vertices has a Hamiltonian cycle if $G$ is of order

$ kp,  1\leq k\leq 47;\quad  kpq,   1\leq k\leq 7;\quad  pqr;  \quad kp\sp 2,  1\leq k\leq 4; $

$kp^3,\,  1\leq k\leq 2; \quad p^k, 1\leq k <\infty,$

where $p$, $q$ and $r$ are distinct primes.
\end{proposition}

\begin{proposition}
	\label{pro2.5}
{\rm \cite[Theorem~1.4]{MR4328721}} With the exception of the Petersen graph, a connected vertex-transitive graph of order $pq$, where $p$ and $q$ are primes, contains a Hamilton cycle.
\end{proposition}

\begin{proposition}
	\label{pro2.6}
{\rm \cite[Theorem~1.1]{10p},\cite[Theorem~1.1]{6p}} Every connected vertex-transitive graph of order $6p$ or $10p$, where $p$ is a prime, contains a Hamilton path. Moreover, with the exception of the truncation of the Petersen graph, every such graph contains a Hamilton cycle.
\end{proposition}

\begin{proposition}
	\label{pro2.7}
{\rm \cite[Theorem~2.4]{A81}} Let $\Gamma$ be a graph that admits a semiregular automorphism $\rho$ of order $t \ge 4$ and $\mathcal{B}$ the set of $\langle \rho \rangle$-orbits on $V(\Gamma)$. Let $\Gamma_1,\Gamma_2,\dots,\Gamma_m$ be the subgraphs of $\Gamma$ induced by the orbits of $\langle \rho \rangle$, and let each $\Gamma_i$ be connected. If $d(\Gamma_j) \ge 3$ for some $j$ and there is a Hamilton path in the quotient graph $\Gamma_{\mathcal{B}}$ with $\Gamma_j$ as one of its endvertices, then $\Gamma$ contains a Hamilton cycle.
\end{proposition}

\begin{proposition}
	\label{pro2.8}
{\rm \cite[Theorem~2]{metacy}} Let $\Gamma$ be a $(m,n)$-metacirculant graph and $G\cong \ZZ_n\rtimes \ZZ_t$ be a metacyclic transitive subgroup of $\Aut(\Gamma)$, where $mn\ge 3$. If $n$ is a prime power, $\Gamma$ is connected, and $\Gamma$ is not the Petersen graph, then $\Gamma$ possesses a Hamilton cycle.
\end{proposition}

\begin{proposition}
	\label{pro2.9}
{\rm \cite[Theorem~2]{Tan90}} Let $\Gamma$ be a connected cubic $(2,n)$-metacirculant graph with automorphisms $\rho$ and $\sigma$ as in the definition, where $n$ is odd. If the subgraph induced by an orbit of $\lg\rho\rg$ contains an edge and $\Gamma$ is not the Petersen graph, then $\Gamma$ contains a Hamilton cycle.
\end{proposition}

\begin{proposition}
	\label{pro2.10}
{\rm \cite[Lemma~2.7]{6p}} Let $\Gamma=X(T,H,D)$ be a coset graph of order $t$ with a Hamilton cycle. Let $C$ be  a Hamilton cycle of $\Gamma$. Let $G =T\times \langle c \rangle $, where $|c|=p$ for a prime coprime to $t$. Then both  $\Gamma_1=X(G,H,Dc\cup Dc^{-1})\cong X(G,H,Dc^k\cup Dc^{-k})$ and $\Gamma_2 = X(G,H,D\cup HcH \cup Hc\sp{-1}H)\cong X(G,H,D\cup Hc^kH \cup Hc\sp{-k}H)$ contain a Hamilton cycle, where \(1\le k\le p-1\). In particular, suppose that $\Gamma_0=X(G,H,D_0)$ is connected and $T$ is primitive on the set $[T:H]$ and $|G:H|=2rs$, where $r$ and $s$ are distinct primes, and $r,s>5$. Then $\Gamma_0$ contains a $G$-minimal graph isomorphic to $\Gamma_1$ or $\Gamma_2$ for a suitable choice of \(D\) and so it has a Hamilton cycle.
\end{proposition}

\begin{lemma}
	\label{le2.11}
Let $\Gamma=X(G,H,D)$ be a $G$-minimal graph that contains a Hamilton cycle. Then every edge in $E(\Gamma)$ is in a Hamilton cycle of $\Gamma$.
\end{lemma}

\demo Let $\Gamma_{0}=X(G,H,D_0)$, where $D_0\subsetneq D$ and $D_0=D_0^{-1}$. Let $D_{g_0}=Hg_0H\cup Hg_0^{-1}H$. We may assume that $D=D_0\cup D_{g_0}$. For any $g_0\in D$, $\Gamma_{0}$ is not connected and so $\Gamma_{0}$ does not contain a Hamilton cycle. That is,  there exists a Hamilton cycle $C$ of $\Gamma$ such that $C$ contains the edge $\{Hg_1,Hg_2\}$, where $g_1g_2^{-1}\in D_{g_0}$. For any edge $e$ in the same $G$-orbit, there exists $g\in G$ taking $\{Hg_1,Hg_2\}$ to $e$. Hence $C^g$ is a Hamilton cycle containing $e$.\qqed

\subsection{Group- and  finite field-theoretic results}
\begin{proposition}
	\label{pro2.12} {\rm \cite[Theorem~1.49]{fsg}}
Every transitive group $G$ of prime degree $p$ is $\ZZ_p\rtimes \ZZ_r$ which is isomorphic to a subgroup of $\AGL(1,p)$, where $r\di (p-1)$;  	$\soc(G)=A_p$;   $\soc(G)=\PSL(2,11)$ of degree 11; $\soc(G)=M_{11}$ or $M_{23}$ of degree $11$ or $23$, respectively; or $\soc(G)=\PSL(d,q)$  of degree $p= (q\sp d-1)/(q-1)$.
\end{proposition}

By Proposition~\ref{pro2.12}, the following result can be extracted from \cite[Proposition 2.5]{MS2}.

\begin{corollary}
	\label{co2.13}
Let $p,q$ be primes and $p> q$ and let $G$ be a quasiprimitive and imprimitive group of degree $pq$. Then $G$ is an almost simple group and $\soc(G)\cong T$, where $T\in \{A_p, \PSL(2,11),M_{11}, M_{23},\PSL(d,r)\}$, where $d$ is a prime, $(d,r-1)=1$ and $p= (r\sp d-1)/(r-1)$.
\end{corollary}

\begin{proposition}
	\label{pro2.14} {\rm \cite{LS85}}
Every primitive group $G$ of degree $2p$, where $p$ is an odd prime, has one of the following socles: $A_{2p}$; $\PSL(2,q)$, where $p=\frac{q+1}2$; $M_{22}$; or $A_5$, where $p\in \{3,5\}$. Except for $A_5$ of degree $10$, all these socles are $2$-transitive.
\end{proposition}

\begin{proposition}
	\label{pro2.15}{\rm \cite{LS85}}
Let $p$ and $q$ be primes with $q<p$. Then any simply primitive group of degree $pq$ is an almost simple group, with socle as listed in {\it Table~\ref{tab1}}.
	\end{proposition}

\begin{table}[htbp]

	\centering
	\caption{ \it Simply primitive groups of degree $pq$ with primes $q<p$}
	\vskip 3mm
	\resizebox{\textwidth}{!}{
		\begin{tabular}{|c|l|c|l|l|}
			\hline
			Row & $T = \operatorname{soc}(G)$ & $(p,q)$ & Action & Remark \\
			\hline
			1 & $A_p$ & $(p,(p-1)/2)$ & pairs & \\
			\hline
			2 & $A_{p+1}$ & $(p,(p+1)/2)$ & pairs & $p\ge5$ \\
			\hline
			3 & $A_7$ & $(7,5)$ & triples & \\
			\hline
			4 & $\operatorname{PSL}(2,7)$ & $(7,3)$ & cosets of $D_{16}$ & $G=\PGL(2,7)$ \\
			\hline
			5 & $\operatorname{PSL}(4,2)$ & $(7,5)$ & $2$-spaces & \\
			\hline
			6 & $\operatorname{PSL}(5,2)$ & $(31,5)$ & $2$-spaces or $3$-spaces & two representations \\
			\hline
			7 & $\mathrm{PSp}(4,2^{2^t})$ & $(2^{2^{t+1}}+1,2^{2^t}+1)$ & $1$-spaces or totally isotropic $2$-spaces & two representations; $t\ge1$ \\
			\hline
			8 & $P\Omega^\pm(2d,2)$ & $(2^{d}\mp 1,2^{d-1}\pm 1)$ & singular $1$-spaces & \\
			\hline
			9 & $\operatorname{PSL}(2,q^2)$ & $((q^2+1)/2,q)$ & cosets of $\operatorname{PGL}(2,q)$ & two representations \\
			\hline
			10 & $\operatorname{PSL}(2,p)$ & $(p,(p\mp 1)/2)$ & cosets of $D_{p\pm 1}$ & \\
			\hline
			11 & $\operatorname{PSL}(2,p)$ & $(19,3)$, $(29,7)$, $(59,29)$ or $(61,31)$ & cosets of $A_5$ & two representations \\
			\hline
			12 & $\operatorname{PSL}(2,23)$ & $(23,11)$ & cosets of $S_4$ & two representations \\
			\hline
			13 & $\operatorname{PSL}(2,11)$ & $(11,5)$ & cosets of $S_4$ & $G=\PGL(2,11)$ \\
			\hline
			14 & $\operatorname{PSL}(2,13)$ & $(13,7)$ & cosets of $A_4$ & \\
			\hline
			15 & $M_{11}$ & $(11,5)$ & pairs & \\
			\hline
			16 & $M_{22}$ & $(11,7)$ & cosets of $\mathbb{Z}_2^4:A_6$ & \\
			\hline
			17 & $M_{23}$ & $(23,11)$ & pairs or cosets of $\mathbb{Z}_2^4:A_7$ & \\
			\hline
		\end{tabular}
	}
	\label{tab1}
\end{table}

The following result can be extracted from \cite[Lemma 4.7]{PWX93}.

\begin{corollary}
	\label{co2.16}
Let $G=\PSL(d,r)$ be a transitive imprimitive permutation group of degree $pq$, where $d\ge 3$ and $p>q>5$ are primes. Let $H$ be a point stabilizer in this action. Suppose that $H<M<G$, where $M$ is the stabilizer of a one-dimensional subspace of $\mathbb{F}_r^d$ and $|G:M|=p$. Then $G$ has exactly $q$ suborbits of length $1$ and a single suborbit of length $q(p-1)$. In particular, $G$ has rank $q+1$.
\end{corollary}

\begin{proposition}
	\label{pro2.17}
{\rm \cite[I. Satz~7.8]{fgt}} If $G$ is a finite group with a normal subgroup $H$ and $P$ is a Sylow-$p$ subgroup of $H$, then $G=N_G(P)H$.
\end{proposition}

\begin{proposition}
	\label{pro2.18}
{\rm \cite[I. Satz~4.5]{fgt}} For a subgroup $H$ of the group $G$, the factor group $N_G(H)/C_G(H)$ is isomorphic to a subgroup of $\Aut(H)$, the group of automorphisms of $H$.
\end{proposition}

\begin{proposition}
	\label{pro2.19} {\rm \cite[I. Satz~17.5]{fgt}}
Let $K$ be an abelian normal subgroup of the group $G$ such that $K\le B\le G$ where $(|K|,|G:B|)=1$. Then $K$ has a complement in $G$ provided $K$ has a complement in $B$.
\end{proposition}

\begin{lemma}
	\label{le2.20}
Let $\Gamma=X(G,H,D)$ be a $G$-minimal graph. Let $\tau\in N_G(H)\backslash H$ and $\tau^{n}\in H$, where $n$ is the smallest positive integer. If $D=\tau^{-j}D\tau^{j}$ and $(n,|Z(G)|)=1$, where $j\in \ZZ_n\backslash \{0\}$, then $\Aut(\Gamma)$ has a transitive subgroup isomorphic to $\ZZ_{\frac{n}{(n,j)}}\times G$.
\end{lemma}

\demo Let $f_j$ be a bijection of $V(\Gamma)$ such that $f_j(Hg)=H\tau^{j}g$. Suppose that $Hg_1\sim Hg_2$. Then $g_1g_2^{-1}\in D$ and $\tau^j g_1\in D\tau^j g_2$, that is, $Hg_1\sim Hg_2$ in $\Gamma$ if and only if $f_j(Hg_1)\sim f_j(Hg_2)$ in $\Gamma$. Thus, $f_j\in \Aut(\Gamma)$. Clearly, $[\lg f_j\rg, G]=1$ and $\lg f_j \rg \cap G=1$, so $\lg f_j \rg \times G$ is a transitive subgroup of $\Aut(\Gamma)$.\qqed

\section{The proof of Theorem~\ref{the:main}}
\subsection{Outline of the proof}
Let $\Gamma=X(G,H,D)$ be a connected non-quasiprimitive graph of order $2pq$ with an intransitive normal subgroup inducing prime-length orbits, where $p$ and $q$ are distinct primes. Suppose $N\triangleleft G$ induces blocks of length $r$, where $r\in \{2,p,q\}$. Let $\B$ be the set of $N$-blocks. By Corollary~\ref{co2.3}, we only need to consider that $N\lessapprox S_r$. Let $P$ be a Sylow-$r$ subgroup of $N$. Then $P\cong \ZZ_r$ and $P$ acts regularly on every $N$-block. Since $N\triangleleft G$ and $P\in \syl_r(N)$, by Proposition~\ref{pro2.17}, we know $G=N_G(P)N$. Thus, $G_0=N_G(P)$ is transitive on $\B$ and is also transitive on $V(\Gamma)$. Let $N_0$ be the kernel of the action of $G_0$ on $\B$. Then $P\le N_0\lessapprox \AGL(1,r)$ and $\overline{G_0}=G_0/N_0$ acts naturally and faithfully on $\B$. When passing to a minimal spanning subgraph, we retain the original vertex action of the transitive group under consideration. All uses of non-quasiprimitivity refer to the original graph. We need to consider the following cases, separately:

\begin{enumerate}
\item $N_0$ induces blocks of length $r$ and $\overline{G_0}$ acts quasiprimitively on $\B$;
\item $N_0$ induces blocks of length $r$ and $\overline{G_0}$ acts non-quasiprimitively on $\B$.
\end{enumerate}

By \cite{KM08} and Proposition~\ref{pro2.6}, we may assume $p>q>5$.

\subsection{$r=p$}

Let $\Gamma=X(G,H,D)$. Under the condition in the outline, in this section we may assume that $G$ contains a normal $p$-subgroup $P\cong \ZZ_p$. Let $\B$ be the set of $P$-blocks and $N$ be the kernel of the action of $G$ on $\B$. Then $G\cong N.G/N$. Set $\overline{G}=G/N$. $G/N$ acts faithfully and transitively on $\B$. Let $H=G_\a$, where $\a\in V(\Gamma)$.
\subsubsection{$\overline{G}$ acts quasiprimitively on $\B$}
In this case, by Corollary~\ref{co2.13} and Proposition~\ref{pro2.14}, we know that $\overline{G}$ is a nonabelian almost simple group and so $\soc(\overline{G})$ is also transitive on $\B$.  We may assume that $\overline{G}=\soc(\overline{G})$.

\begin{lemma}
	\label{le3.1}
$C_G(P)$ is transitive on $V(\Gamma)$.
\end{lemma}

\demo  By Proposition~\ref{pro2.18}, we get that $G/C_G(P)\lessapprox \Aut(P)\cong \ZZ_{p-1}$, which implies that $C_G(P)\neq P$. Since $P\le C_G(P)\triangleleft G$, if $C_G(P)$ is not transitive on $V(\Gamma)$, then $C_G(P)N$ is also not transitive. Since $\overline{G}$ acts quasiprimitively on $\B$, we get that $C_G(P)$ acts trivially on $\B$ and $C_G(P)=P$, a contradiction.\qqed
\begin{lemma}
	\label{le3.2}
$G_0:=C_G(P)\cong P\times\overline{G}$.
\end{lemma}

\demo Since $C_G(P)$ is transitive on $V(\Gamma)$, we know that $1\neq C_G(P)N/N\triangleleft \overline{G}$. Since $\overline{G}$ is a simple group, we get that $C_G(P)N/N=\overline{G}$, $C_G(P)N=G$. $C_G(P)/P\cong C_G(P)N/N\cong \overline{G}$. Let $B_0\in \B$ such that $\a\in B_0$. Then $(G_0)_{B_0}=P\times H_0$, where $H_0=G_0\cap H$. Since $(|P|,|G_0:PH_0|)=(p,2q)=1$, by Proposition~\ref{pro2.19}, we get that $P$ has a complement in $G_0$. Thus, $G_0\cong P\times \overline{G}$. We identify $\overline{G}$ with the complement of $P$ in $G_0$. Since the original graph is non-quasiprimitive, $H_0\le\overline{G}$; otherwise, $\overline{G}$ would act transitively on its vertices.\qqed

\begin{theorem}
Let $\Gamma_0$ be a $G_0$-minimal spanning subgraph of $\Gamma$. Then $\Gamma_0$ has a Hamilton cycle.
\end{theorem}

\demo Write $\Gamma_0=X(G_0,H_0,D_0)$, where $H_0=G_0\cap H$. For simplicity, in this proof we write $\Gamma,H,D$ for $\Gamma_0,H_0,D_0$, respectively. Suppose that $\overline{G}$ acts primitively on $\B$, as listed in Proposition~\ref{pro2.14}. Then by Propositions~\ref{pro2.5} and \ref{pro2.10}, we get that $\Gamma$ has a Hamilton cycle. Thus, we only need to consider the case where $\overline{G}$ acts quasiprimitively and imprimitively on $\B$. By Corollary~\ref{co2.13}, we know that the other almost simple groups in Corollary~\ref{co2.13} are excluded by $q>5$, and so either $\overline{G}\cong \PSL(d,q_0)$, where $q=(q_0^d-1)/(q_0-1)$; or $\overline{G}\cong M_{11}$. In the first case, let $M$ be a maximal subgroup of $\overline{G}$ containing $H$ and $g_0\in \overline{G}\backslash M$. Then $|M:H|=2$ and $\overline{G}=M\cup Mg_0M=M\cup Mg_0H$. Let $\ell\in M\backslash H$. Then $\ell^2\in H$ and $\lg H,\ell\rg=M$. In particular, $\overline{G}=H\cup H\ell\cup Hg_0H$. Thus, $H\ell g_0H=Hg_0\ell H=Hg_0H$. Set $\lg c\rg =P$. Now, we give all forms of the double cosets of $G_0$ with respect to $H$:

\begin{enumerate}
\item $Hc^iH$, where $i\in \ZZ_p$;
\item $Hc^i g_0H$, where $i\in \ZZ_p$;
\item $Hc^i\ell^jH$, where $i\in \ZZ_p$ and $j\in \ZZ_2$.
\end{enumerate}

Since $\lg H,c,\ell\rg = P\times M\neq G_0$, we know that $D$ must contain a double coset of case 2. Suppose $Hc^ig_0H\subset D$, where $i\neq 0$. Then $\lg H,c^ig_0\rg= G_0$. In this case, $D=c^iH g_0H\cup c^{-i}Hg_0H$. Let $\Gamma_1=X(\overline{G},H,H g_0H)$, a connected vertex-transitive graph, containing a Hamilton cycle. By Proposition~\ref{pro2.10}, $\Gamma$ contains a Hamilton cycle. Thus, we may assume that $D$ does not contain this coset. Suppose $Hg_0H\subset D$. Since $\lg H,\ell, g_0\rg=\overline{G}$, we may assume that $D$ contains a double coset which is $Hc\ell^kH$, where $k\in \ZZ_2$. Since $\lg H, g_0, c\ell^k\rg= G_0$, we get  $D=Hc\ell^kH\cup Hc^{-1}\ell^kH\cup Hg_0H$. Let $\Gamma_2=X(\overline{G},H,H g_0H)$. In particular, $D^\ell= D$. By Lemma~\ref{le2.20}, there exists a subgroup $G_1$ of $\Aut(\Gamma)$ isomorphic to $\ZZ_2\times {G_0}\cong \ZZ_{2p}\times \overline{G}$ such that $G_1$ is transitive on $\Gamma$. The extra $\ZZ_2$ interchanges the two members of each of the $q$ blocks of size $2$ on $\B$. Since an element of order $q$ in $\overline{G}$ cyclically permutes these blocks, $G_1$ contains a subgroup isomorphic to $\ZZ_{2pq}$ acting regularly on $V(\Gamma)$. By Proposition~\ref{pro2.4}, $\Gamma$ contains a Hamilton cycle.

Now we suppose that $G_0= P\times M_{11}$, where $P=\lg c\rg\cong \ZZ_p$. Then $H\cong A_6$ and we may assume $H<M_{10}<M_{11}$. Let $\ell\in M_{10}\backslash H$, and let $g_0\in M_{11}\backslash M_{10}$. Then $M_{11}=H\cup H\ell H\cup H g_0H$. We give all forms of the double cosets of $G_0$ with respect to $H$:

\begin{enumerate}
	\item $Hc^iH$, where $i\in \ZZ_p$;
	\item $Hc^i g_0H$, where $i\in \ZZ_p$;
	\item $Hc^i\ell^jH$, where $i\in \ZZ_p$ and $j\in \ZZ_2$.
\end{enumerate}

Using the same argument as above, we get that $\Gamma$ contains a Hamilton cycle.\qqed

\subsubsection{$\overline{G}$ acts non-quasiprimitively on $\B$}
In this case, let $\overline{T}$ be a maximal intransitive normal subgroup of $\overline{G}$. Then $\overline{T}$ induces blocks of length $t$, where $t\in \{2,q\}$. Thus, we consider two cases $t=2$ and $t=q$, respectively. Suppose $N\cong \ZZ_p\rtimes \ZZ_s\lessapprox \AGL(1,p)$. Firstly, we consider a special case in Lemma~\ref{le3.4}.

\begin{lemma}
	\label{le3.4}
If $G= \ZZ_{pq}\rtimes \ZZ_{2^i}$, then $\Gamma$ contains a Hamilton cycle.
\end{lemma}
\demo In this case, $H=\ZZ_{2^{i-1}}$. If $i=1$, then $\Gamma$ is a Cayley graph. By Proposition~\ref{pro2.4}, $\Gamma$ contains a Hamilton cycle. If $\ZZ_p\le Z(G)$, then $G\cong \ZZ_q\rtimes \ZZ_{2^ip}$. By Proposition~\ref{pro2.8}, $\Gamma$ contains a Hamilton cycle. Similarly, the existence of a Hamilton cycle in \(\Gamma\) still holds when $\ZZ_q\le Z(G)$. If $\ZZ_2\le Z(G)$, then $G$ contains a subgroup isomorphic to $\ZZ_{2pq}$. By Proposition~\ref{pro2.4}, $\Gamma$ contains a Hamilton cycle. For the rest of the proof, we suppose that $Z(G)=1$ and $i\ge2$.

Let $C=\lg x\rg\cong \ZZ_{pq}$ be the cyclic subgroup of $G$ of order $pq$. Then $P=\lg x^q\rg $ and $Q=\lg x^p\rg$ are the Sylow-$p$ subgroup and Sylow-$q$ subgroup, respectively. Let $M=\lg \ell\rg \le G$ and $H=\lg \ell^2\rg$. Then $PQ$ induces blocks of length $pq$. The kernel of the action of $G$ on this block system is $T=PQH=\ZZ_{pq}\rtimes\ZZ_{2^{i-1}}$. Since $P\triangleleft G$ and $Q\triangleleft G$, we know that these two groups induce blocks of length $p$ and $q$, respectively. Let $\Q$ be the set of $Q$-blocks. Let $\a\in B_1\in \B$ and $\a \in B_2\in \Q$. Let $\P_1=\{B_1^g:g\in T\}$ and $\P_2=\{B_1^{\ell g}:g\in T\}$. Let $\Q_1=\{B_2^g:g\in T\}$ and $\Q_2=\{B_2^{\ell g}:g\in T\}$. Let $N_{\P_1}$ and $N_{\Q_1}$ be the kernels of the actions of $T$ on $\P_1$ and $\Q_1$, respectively.  Let $\P=\B$, and write $N_\P=N_{\P_1}$ and $N_\Q=N_{\Q_1}$ for the kernels of the actions of $G$ on $\P$ and $\Q$, respectively. To prove this lemma, we consider the following three parts.

 \vskip 3mm
{\it Part 1: $N_{\P_1}=P$ and $N_{\Q_1}=Q$.}
 \vskip 3mm
In this case, we claim that $C_T(P)=C_T(Q)=PQ$. Clearly, $PQ \triangleleft C_T(P)$. If $PQ\neq C_T(P)$, then there exists $g_0\in C_T(P)$ with $o(g_0)=2^j$, where $1\le j\le i-1$. Thus, $C_T(P)\cong P\times (Q\rtimes \lg g_0\rg )$. Since $Q\rtimes \lg g_0\rg \char C_T(P) \triangleleft T$, we get that $Q\rtimes \lg g_0\rg\triangleleft T$. Since $|Q\rtimes \lg g_0\rg|=2^jq$ and it is not transitive on $\{Hg:g\in T\}$, we know that $Q\rtimes \lg g_0\rg\le N_{\Q_1}$, a contradiction. Using the same argument as above, we get $C_T(P)=C_T(Q)=PQ$.  Let $\ell_0=\ell^{2^{i-1}}$. Then $(x^{p})^{\ell_0}=x^{-p}$ and $(x^{q})^{\ell_0}=x^{-q}$ and so $x^{\ell_0}=x^{-1}$. Now, we give all forms of the double cosets of $G$ with respect to $H$:

\begin{enumerate}
\item $Hx^{kp}H$, where $k\in \ZZ_q\backslash \{0\}$, is self-paired and the length of it is $2^{i-1}$;
\item $Hx^{kq}H$, where $k\in \ZZ_p\backslash \{0\}$, is self-paired and the length of it is $2^{i-1}$;
\item $Hx^{k_1p+k_2q}H$, where $k_1\in \ZZ_q$, $k_2\in \ZZ_p$ and $k_1k_2\neq 0$, is self-paired and the length of it is $2^{i-1}$;
\item $H\ell x^{kp}H$, where $k\in \ZZ_q\backslash \{0\}$, is not self-paired and the length of it is $2^{i-1}$;
\item $H\ell x^{kq}H$, where $k\in \ZZ_p\backslash \{0\}$, is not self-paired and the length of it is $2^{i-1}$;
\item $H\ell x^{k_1p+k_2q}H$, where $k_1\in \ZZ_q$, $k_2\in \ZZ_p$ and $k_1k_2\neq 0$, is not self-paired and the length of it is $2^{i-1}$;
\item $H\ell H$ is self-paired and the length of it is 1.
\end{enumerate}
We consider all forms of the $G$-minimal graph.

Firstly, since $\lg H, \ell x^{k_1p+k_2q}\rg=G$, where $k_1\in \ZZ_q$, $k_2\in \ZZ_p$ and $k_1k_2\neq 0$, we get that $\Gamma=X(G,H,H\ell x^{k_1p+k_2q}H\cup H x^{-(k_1p+k_2q)} \ell H)$ is a $G$-minimal graph. Let $x_0=x^{k_1p+k_2q}$. Then $\lg x_0\rg =PQ$. Since $H\ell x_0 \ell^{2^{i-1}}=H\ell x_0^{-1}\in N(H)$, we get a path of length 2 in $\Gamma$, which is

 $$H\ell x_0,H,H\ell x_0^{-1}.$$

Acting on this path by $\lg x_0^{-2}\rg$, we get the Hamilton cycle of $\Gamma$, which is

$$H\ell x_0,H,H\ell x_0^{-1},Hx_0^{-2},\ldots, Hx_0^{4},H\ell x_0^3,Hx_0^2,H\ell x_0.$$

Secondly, we consider that $H\ell x^{kp}H\subset D$. We may assume that $k=1$. Since $\Gamma$ is a $G$-minimal graph, we know that $D$ cannot contain any double coset of the form appearing in Case 6. Since $\lg H,\ell, x^{p}, \ell x^{p}\rg\neq G$, we know that $D$ contains a double coset of the form appearing in Case 2, 3 or 5. For the three cases, $\Gamma_\Q\cong X(\overline{G},\overline{H},\overline{D})$, where $\overline{G}=G/N_\Q$. In particular, $\Gamma_\Q$ is a $\overline{G}$-minimal graph of order $2p$ if $D$ contains a double coset of Case 2 or 3. Since $q>5$, $\Gamma_\Q$ contains a Hamilton cycle. Let $B_1=\{Hg:g\in HQ\}$ and $B_2=\{Hg:g\in H\ell Q\}$. Since $B_2\cap N(H)$ contains a subset $\{H\ell x^{p},{H\ell x^{-p}}\}$ (a set not an edge), we know that $d(X[B_1,B_2])\ge 2$. By Proposition~\ref{pro2.2} and Lemma~\ref{le2.11}, $\Gamma$ contains a Hamilton cycle. Using the same argument as above, we know that if $D$ contains a double coset of the form appearing in Case 2 or 3, then $\Gamma$ contains a Hamilton cycle. If $D$ contains a double coset of Case 5, then $\Gamma_\Q$ contains a Hamilton cycle, which is
$$
\overline{H},\ \overline{H}\ell,\
\overline{H}x^{-q},\ \overline{H}\ell x^{-q},\
\overline{H}x^{-2q},\ \overline{H}\ell x^{-2q},\ \ldots,\
\overline{H}x^{-(p-1)q},\
\overline{H}\ell x^{-(p-1)q},\ \overline{H}.
$$
Since $d(X[\overline{H},\overline{H}\ell])\ge 2$, by Proposition~\ref{pro2.2}, $\Gamma$ contains a Hamilton cycle.

Thirdly, we consider that $Hx^{k_1p+k_2q}H\subset D$, where $k_1\in \ZZ_q$, $k_2\in \ZZ_p$ and $k_1k_2\neq 0$. Let $x_0=x^{k_1p+k_2q}$. Then $\lg x_0\rg=\lg x\rg$. Since $\lg H, x^p,x^q, x\rg\neq G$, we know that $D=Hx_0H\cup H\ell H$. Suppose $i=2$. Then $G=\lg x,\ell:x^{pq}=\ell^4=1, x^\ell=x^t, t^2=-1\pmod {pq}\rg$. Thus, $x^{\ell^2}=x^{-1}$. Since $\lg H, x,x^p,x^q\rg\neq G$, we know that $\Gamma\cong X(G,H,HxH\cup H\ell H)$, which is a cubic $(2,pq)$-metacirculant graph. By Proposition~\ref{pro2.9}, $\Gamma$ contains a Hamilton cycle. Suppose $i\ge 3$. Let $\B$ be the set of $T$-blocks. Then $|\B|=2$. Let $T_1$ and $T_2$ be the two distinct $T$-block graphs. Then $d(T_1)=d(T_2)\ge 4$ and $\Gamma_\B\cong K_2$. Clearly, $T_1$ and $T_2$ are connected Cayley graphs. By Proposition~\ref{pro2.7}, $\Gamma$ contains a Hamilton cycle.

Finally suppose that $D=Hx^pH\cup Hx^qH\cup H\ell H$. Let $\B$ be the set of $T$-blocks and $T_1$ and $T_2$ be the two distinct $T$-blocks. Then using the same argument as above, $\Gamma$ contains a Hamilton cycle.

 \vskip 3mm
{\it Part 2: $N_{\P_1}\neq P$ and $N_{\Q_1}=Q$.}
\vskip 3mm
In this case, we claim that $N_{\P_1}\times Q= C_T(Q)$. Suppose $N_{\P_1}=P\rtimes \lg \ell^{2^t}\rg$, where $1\le t\le i-1$. Since $N_{\P_1}\cap Q=1$, we know that $N_{\P_1}\times Q\le C_T(Q)$. Since $Q\le Z(C_T(Q))$ and $Q\in \syl_q(C_T(Q))$, we know that $Q$ has a unique normal complement in $C_T(Q)$, denoted by $C$. It remains to show $C=N_{\P_1}$. Let $H\ell x^{jp}P\in\P_2$ be any block in $\P_2$. Then for any $g\in C$, $H\ell x^{jp}Pg=\ell H gPx^{jp}=\ell H Px^{jp}=H\ell x^{jp}P$. Thus, $C\le N_{\P_1}$ and $N_{\P_1}\times Q= C_T(Q)$.

Let $\ell_0=\ell^{2^{i-1}}$ and $\ell_j=\ell^{2^j}$, where $2\le j\le i$. Then $(x^{q})^{\ell_{0}}=x^{-q}$ and \((x^p)^{\ell^{2^{t-1}}}=x^{-p}\). Now, we give all forms of the double cosets of $G$ with respect to $H$:

\begin{enumerate}
\item $Hx^{kp}H$, where $k\in \ZZ_q\backslash \{0\}$, is self-paired if $t\ge 2$ and is not self-paired if $t=1$, and the length of it is $2^{t-1}$;
\item $Hx^{kq}H$, where $k\in \ZZ_p\backslash \{0\}$, is self-paired and the length of it is $2^{i-1}$;
\item $Hx^{k_1p+k_2q}H$, where $k_1\in \ZZ_q$, $k_2\in \ZZ_p$ and $k_1k_2\neq 0$, is not self-paired and the length of it is $2^{i-1}$;
\item $H\ell x^{kp}H$, where $k\in \ZZ_q\backslash \{0\}$, is self-paired if $t=1$ and not self-paired if $t\ge 2$, and the length of it is $2^{t-1}$;
\item $H\ell x^{kq}H$, where $k\in \ZZ_p\backslash \{0\}$, is not self-paired and the length of it is $2^{i-1}$;
\item $H\ell x^{k_1p+k_2q}H$, where $k_1\in \ZZ_q$, $k_2\in \ZZ_p$ and $k_1k_2\neq 0$, is not self-paired and the length of it is $2^{i-1}$;
\item $H\ell H$ is self-paired and the length of it is 1.
\end{enumerate}

Firstly, we consider that $H\ell x^{k_1p+k_2q}H\subset D$, where $k_1\in \ZZ_q$, $k_2\in \ZZ_p$ and $k_1k_2\neq 0$. Similarly, let $x_0= x^{k_1p+k_2q}$. If $t\ge 2$, then $\lg H,\ell x_0\rg= G$. Thus, $\Gamma=X(G,H,H\ell x_0H\cup H x_0^{-1} \ell H)$. Using the same argument as in Part 1, we get that $\Gamma$ contains a Hamilton cycle. Suppose $t=1$. Then $\lg H,\ell x_0\rg\neq G$. Since $\Gamma$ is connected, there exists $g_0\in G\setminus\lg H,\ell x_0\rg$ such that $D$ contains a double coset $Hg_0H$. Since $q>5$, we know that $\Gamma_\P\cong X(\overline{G},\overline{H},\overline{H\ell x_0H}\cup\overline{Hx_0^{-1}\ell H}\cup\overline{Hg_0 H}\cup \overline{Hg_0^{-1}H})$, which is a $\overline{G}$-minimal graph. Using the same argument as Part 1, we get that $\Gamma$ contains a Hamilton cycle.

In fact, if $D$ contains a double coset of the form appearing in Case 4 or 5, then using the same argument as in Part 1, we get that $\Gamma$ contains a Hamilton cycle.

Secondly, if $D$ contains a double coset of the form appearing in Case 3, using the same argument as in Part 1, we know that the only untreated case is $D=Hx^{k_1p+k_2q}H\cup Hx^{-k_1p-k_2q}H\cup H\ell H$. Since $q>5$, we know that $\Gamma_\P\cong X(\overline{G},\overline{H},\overline{Hx^{k_1p}H}\cup\overline{Hx^{-k_1p}H}\cup\overline{H\ell H})$ contains a Hamilton cycle. Using the same argument as Part 1, we get that $\Gamma$ contains a Hamilton cycle.

Finally, if $D$ contains a double coset of the form appearing in Case 1, then using the same argument as Part 1, we get that $\Gamma$ contains a Hamilton cycle.

Clearly, since we have not used the condition \(p>q\), we know that the case $N_\P= P$ and $N_\Q\neq Q$ is analogous to Part 2.

  \vskip 3mm
{\it Part 3: $N_{\P_1}\neq P$ and $N_{\Q_1}\neq Q$.}
 \vskip 3mm
By Part 2, we know that $\ell_0\in C_T(P)\cap C_T(Q)$. Thus, $\lg x,\ell_0\rg\cong \ZZ_{2pq}$. However, $\ell_0\in Z(G)$, a contradiction with $Z(G)=1$.

By Parts 1-3, we know that if $G\cong \ZZ_{pq}\rtimes \ZZ_{2^i}$, where $p>q>5$, then $\Gamma$ contains a Hamilton cycle. \qqed
\begin{lemma}
	\label{le3.5}
If $t=2$, then $\Gamma$ contains a Hamilton cycle.
\end{lemma}
\demo Let $\overline{M}$ be the Sylow-2 subgroup of $\overline{T}$. Then by Proposition~\ref{pro2.17}, $N_{\overline{G}}(\overline{M})$ is transitive on $\B$. Let $\overline{G_0}=N_{\overline{G}}(\overline{M})$. Then $G$ contains a subgroup $G_0$ isomorphic to $N.\overline{G_0}$ and this subgroup is transitive on $V(\Gamma)$. Since $\overline{M}\triangleleft \overline{G_0}$ and $\overline{M}$ induces blocks of length 2, we get that the exponent of $\overline{M}$ is 2. Thus, $\overline{M}\cong \ZZ_2^n$ and $\overline{G_0}/\overline{M}\lessapprox S_q$, where $n\le q$. Let $\overline{Q}$ be a Sylow-$q$ subgroup of $\overline{G_0}$. Then $\overline{M.Q}$ acts transitively on $\B$. Since $\overline{M}$ is a normal Sylow-$2$ subgroup of $\overline{M.Q}$, we know $\overline{M.Q}\cong \overline{M}\rtimes \overline{Q}$. Thus, $G$ contains a subgroup $G_1$ isomorphic to $N.(\overline{M}\rtimes \overline{Q})\cong N.(\ZZ_2^i\rtimes \ZZ_q)$, where $1\le i\le n$. If $Z(\overline{M}\rtimes \overline{Q})\neq 1$, then there exists an element $\r\in \overline{M}\rtimes \overline{Q}$ such that $o(\r)=2q$. Since $\overline{G}$ acts faithfully on $\B$, we know that $N.\lg \r\rg$ acts transitively on $V(\Gamma)$. Let $G_\rho=N.\lg \rho\rg$. If $C_{G_\rho}(P)=P$, then by Proposition~\ref{pro2.18}, we get $G_\rho\cong \ZZ_p\rtimes \ZZ_{2qs}$. By Proposition~\ref{pro2.8}, $\Gamma$ contains a Hamilton cycle. If $C_{G_\rho}(P)\cong \ZZ_{2p}$, then $Z(G_\rho) \cong \ZZ_2$ and $G_\rho\cong \ZZ_{2p}.\ZZ_{qs}$. Thus, $G_\rho$ contains a subgroup isomorphic to $\ZZ_{2p}.\ZZ_{q^i}\cong \ZZ_p\rtimes \ZZ_{2q^i}$, which acts transitively on $V(\Gamma)$, where $q^i\mid qs$. By Proposition~\ref{pro2.8}, $\Gamma$ contains a Hamilton cycle. Suppose $C_{G_\rho}(P)\cong \ZZ_{pq}$. Let $M_0$ be the Sylow-$2$ subgroup of ${G_\rho}$. Then $C_{G_\rho}(P)\rtimes M_0\cong \ZZ_{pq}\rtimes \ZZ_{2^j}$ acts transitively on $\Gamma$, where $2^j\mid 2s$. By Lemma~\ref{le3.4}, $\Gamma$ contains a Hamilton cycle. If $C_{G_\rho}(P)\cong \ZZ_{2pq}$, then $\Gamma$ is a Cayley graph. By Proposition~\ref{pro2.4}, $\Gamma$ contains a Hamilton cycle.

Suppose that $Z( \overline{M}\rtimes \overline{Q})=1$. Then $i\ge 2$ and $2\mid |C_{G_1}(P)|$. Let $M_0$ be the Sylow-2 subgroup of $C_{G_1}(P)$. Then $M_0\not\le N$. Thus, $C_{G_1}(P)$ either induces blocks of length $2p$, or is transitive on $V(\Gamma)$. We have $C_{G_1}(P)/P\cong C_{G_1}(P)N/N\le G_1/N\cong \ZZ_2^i\rtimes \ZZ_q$. If $q~|~|C_{G_1}(P)|$, then $C_{G_1}(P)$ acts transitively on $V(\Gamma)$. Since $Z( \overline{M}\rtimes \overline{Q})=1$, we know that $C_{G_1}(P)\cong  P\times (\overline{M}\rtimes \overline{Q})\cong \overline{M}\rtimes P\overline{Q}\cong \ZZ_2^i\rtimes \ZZ_{pq}$. Since $i\ge 2$, by Corollary~\ref{co2.3}, $\Gamma$ contains a Hamilton cycle. If $q\nmid |C_{G_1}(P)|$, then $M_0 \char C_{G_1}(P)\triangleleft G_1$ and so $M_0\triangleleft G_1$. If $M_0\cong \ZZ_2$, then $M_0\le Z(G_1)$ and $1\neq M_0N/N\le Z( \overline{M}\rtimes \overline{Q})$, a contradiction. If $M_0\not\cong \ZZ_2$, then $M_0$ induces blocks of length 2. By Corollary~\ref{co2.3}, $\Gamma$ contains a Hamilton cycle. \qqed

\begin{lemma}
	\label{le3.6}
If $t=q$, then $\Gamma$ contains a Hamilton cycle.
\end{lemma}

\demo Let $\P$ be the set of $P$-blocks, so that $|\P|=2q$, and let $\overline{Q}$ be a Sylow-$q$ subgroup of $\overline{T}$. By Proposition~\ref{pro2.17}, $N_{\overline{G}}(\overline{Q})$ is transitive on $\P$. Since $\overline{Q}\lessapprox\ZZ_q^2$ has two orbits of length $q$, we may choose $\ell\in G$ with $o(\ell)=2^n$ such that $\ell N\in N_{\overline{G}}(\overline{Q})$ interchanges these two orbits. Let $G_0$ be the subgroup of $G$ containing $N$ such that
\[
\overline{G_0}=G_0/N=\lg\overline{Q},\ell N\rg.
\]
Then $G_0$ is transitive on $V(\Gamma)$.

Since $C_{G_0}(P)\cap N=P$ and
\[
C_{G_0}(P)/P\cong C_{G_0}(P)N/N
\le\overline{G_0},
\]
we have $(p,|C_{G_0}(P):P|)=1$. By Proposition~\ref{pro2.19}, $P$ has a complement $C$ in $C_{G_0}(P)$. Since $P\le Z(C_{G_0}(P))$, we get $C_{G_0}(P)=P\times C$, and this complement is unique. Thus $C\char C_{G_0}(P)\triangleleft G_0$, and hence $C\triangleleft G_0$ and $C\cap N=1$.

If $C=1$, then by Proposition~\ref{pro2.18}, $G_0/P\lessapprox\Aut(P)$ is cyclic. By Proposition~\ref{pro2.19}, $G_0\cong\ZZ_p\rtimes\ZZ_s$, where $s=|G_0:P|$. Proposition~\ref{pro2.8} gives a Hamilton cycle in $\Gamma$.

Suppose that $q\mid|C|$, and let $M_0$ be a Sylow-$q$ subgroup of $C$. Since $C$ embeds in $\overline{G_0}$ and $\overline{Q}\triangleleft\overline{G_0}$, $M_0$ is unique and hence normal in $G_0$. If $M_0\cong\ZZ_q^2$, then $M_0$ induces blocks of length $q$ and acts unfaithfully on each block. The corresponding quotient graph has order $2p$. By Proposition~\ref{pro2.5} and Corollary~\ref{co2.3}, $\Gamma$ contains a Hamilton cycle. If $M_0\cong\ZZ_q$, then $PM_0\cong\ZZ_{pq}\triangleleft G_0$, and $\ell$ interchanges its two orbits. Thus $\lg P,M_0,\ell\rg\cong\ZZ_{pq}\rtimes\ZZ_{2^n}$ is transitive on $V(\Gamma)$. By Lemma~\ref{le3.4}, $\Gamma$ contains a Hamilton cycle.

Finally, suppose that $C\neq1$ and $q\nmid|C|$. By Proposition~\ref{pro2.18}, $G_0/C_{G_0}(P)\lessapprox\Aut(P)$ is cyclic. Thus, for every $g\in\overline{Q}$,
\[
g^{\ell N}g^{-1}\in
\overline{Q}\cap\bigl(C_{G_0}(P)N/N\bigr)=1.
\]
Therefore $\overline{G_0}$ is abelian. Since it acts faithfully and transitively on $\P$, we get $\overline{G_0}\cong\ZZ_{2q}$. Consequently, $C\cong\ZZ_2$ and $C\le Z(G_0)$. Let $Q$ be a Sylow-$q$ subgroup of $G_0$. Since $q\nmid|C_{G_0}(P)|$, $Q$ embeds in the cyclic group $G_0/C_{G_0}(P)$, and hence $Q\cong\ZZ_{q^i}$ for some $i\ge1$. The images of $Q$ and $C$ modulo $N$ generate $\overline{G_0}$, so $PQC\cong\ZZ_p\rtimes\ZZ_{2q^i}$ is transitive on $V(\Gamma)$. By Proposition~\ref{pro2.8}, $\Gamma$ contains a Hamilton cycle. \qqed

\begin{theorem}
	\label{th:3.7}
If $\overline{G}$ acts non-quasiprimitively on $\B$, then $\Gamma$ has a Hamilton cycle.
\end{theorem}
\demo By Lemmas~\ref{le3.4}-\ref{le3.6}, the theorem is proved.\qqed

\begin{remark}
The case $r=q$ is obtained from Section~3.2 by interchanging $p$ and $q$. Indeed, no argument in the $r=p$ case uses the assumption $p>q$: the quotient $\overline{G}$ acts on $2q$ points (respectively $2p$ points when $r=q$), and every group-theoretic result invoked (Propositions~\ref{pro2.12}--\ref{pro2.15}, Corollary~\ref{co2.13}) is symmetric in the two primes.  Thus the entire proof in Section~3.2 applies verbatim with $p$ and $q$ swapped, and the case $r=q$ is settled.
\end{remark}

\subsection{$r=2$}

Let $\Gamma=X(G,H,D)$. Under the condition in the outline, in this section we may assume that $G$ contains a normal subgroup $N\cong \ZZ_2$. Let $\B$ be the set of $N$-blocks. Then $G\cong N.G/N$. Set $\overline{G}=G/N$. $G/N$ acts faithfully and transitively on $\B$. In particular, $N\le Z(G)$. Let $H=G_\a$, where $\a\in V(\Gamma)$.
\subsubsection{$\overline{G}$ acts quasiprimitively on $\B$}

In this case, by Corollary~\ref{co2.13} and Proposition~\ref{pro2.14}, we know that $\overline{G}$ is also a nonabelian almost simple group and so $\soc(\overline{G})$ is also transitive on $\B$.  We may assume that $\overline{G}=\soc(\overline{G})$. In particular, Lemmas~\ref{le3.1}-\ref{le3.2} also hold when $p=2$. Thus, $G=\ZZ_2.\overline{G}= N\times \overline{G}= \ZZ_2\times \overline{G}$. Set $N=\lg \tau \rg$.

\begin{theorem}
If $\overline{G}$ acts quasiprimitively on $\B$, then $\Gamma$ has a Hamilton cycle.
\end{theorem}

\demo Suppose that $\overline{G}$ acts primitively on $\B$, as listed in Table~\ref{tab1} of Proposition~\ref{pro2.15}. Then by Propositions~\ref{pro2.5} and \ref{pro2.10}, we get that $\Gamma$ has a Hamilton cycle. Thus, we may assume that $\overline{G}$ acts quasiprimitively and imprimitively on $\B$. By Corollary~\ref{co2.13}, we get $\overline{G}\cong\SL(d,r)$. Let $\theta$ be a primitive element of $\FF_r$ so that $\lg \theta \rg=\FF_r^*$.

Firstly, suppose that $d=2$. Then $r=2^{2^s}$ and $r+1$ is a Fermat prime. Without loss of generality, we may assume that $\overline{G}=\SL(2,2^{2^s})$. Let $n=2^s$. Set

$$
u_i=\begin{pmatrix}
	1 & \theta \sp{i}\\
	0 & 1
\end{pmatrix},
u_i \sp{'}=	\begin{pmatrix}
	1 & 0\\
	\theta \sp{i} & 1
	
\end{pmatrix},
\ell_i=	\begin{pmatrix}
	\theta \sp{i} & 0\\
	0 & \theta \sp{-i}
\end{pmatrix},
t_i=	\begin{pmatrix}
	0 & \theta \sp{i}\\
	\theta \sp{-i} & 0
\end{pmatrix},
$$
where $\theta$ is a primitive element of $\FF_{2^n}$. Then $H=\lg u_0,u_1,\ell_q\rg\cong \ZZ_2^{n}\rtimes \ZZ_\frac{2^n-1}{q}$. Let $M= \lg u_0,u_1,\ell\rg\cong \ZZ_2^{n}\rtimes \ZZ_{2^n-1}$, where $\ell=\ell_1$. Then $M$ is the unique maximal subgroup of $\overline{G}$ containing $H$ and $\overline{G}=M\cup Mt_0M=Ht_0M\cup HM$. We give all forms of the double cosets of $G$ with respect to $H$:

\begin{enumerate}[(i)]
	\itemsep=0pt

\item $H\tau^i\ell^j t_0H$, where $i\in \ZZ_2$ and $j\in \ZZ_q$, $2q$ suborbits of length $p-1$, all of them are self-paired.

\item $H\tau^i \ell^jH$, where $(i,j)\neq (0,0)$, $2q-1$ suborbits of length $1$, all of them are not self-paired except for $H\tau H$.

\end{enumerate}

For any $i\in \FF_q$, $\lg H,\tau \ell^it_0\rg=G$. If $H\tau\ell^i t_0H = D$, then by Proposition~\ref{pro2.10}, $\Gamma$ contains a Hamilton cycle. Let $D_i=H\ell^i t_0H= H t_0\ell^{-i}H$ and $D_j'=H\tau \ell^jH\cup H\tau\ell^{-j}H$. Since $\lg H,\ell^it_0\rg =\overline{G}$ and $\lg H, \tau\ell^j\rg =N\times M$, where $j\neq 0$, we know that if $D\neq H\tau\ell^i t_0H$, then $D=D_i\cup D_j'$, where $i,j\in \ZZ_q$. If $j=0$, then by Proposition~\ref{pro2.10}, $\Gamma$ contains a Hamilton cycle. Now we suppose that $j\neq0$. Let $B_g=\{Hx:x\in Mg\}$ and $\B'=\{B_g:g\in G\}$. Set $\B_1'=\{B_g:g\in \overline{G}\}$ and $\B_2'=\{ B_{\tau g}:g\in \overline{G} \}$. Then $\B'$ is a  system of blocks. Consider the quotient graph $\Gamma_{\B'}$, whose vertex set is $\B'$. $B_{g_1}\sim B_{g_2}$ in $\Gamma_{\B'}$ if and only if $X[B_{g_1},B_{g_{2}}]$ contains a matching. We claim that the induced subgraphs $\Gamma_{\B_1'}$ and $\Gamma_{\B_2'}$ are complete graphs of order $p$. Let $Hg_1$ and $Hg_2$ be two distinct vertices of $\Gamma$, where $g_2g_1^{-1}\in D_i$. Then $\{Hg:g\in D_i\ell^mg_1\}=\{Hg:g\in H\ell^{i-m}t_0Hg_1\}\subset N(H\ell^mg_1)$. Thus, if $Hg_1\sim Hg_2$, then $H\ell^mg_1\sim H\ell^{-m}g_2$. $H\ell^{-m_1}g_2=H\ell^{-m_2}g_2$ if and only if $\ell^{m_1-m_2}\in H$, that is, $H\ell^{m_1}=H\ell^{m_2}$. Thus, $X[B_{g_1},B_{g_2}]$ contains a matching. Since $\overline{G}$ acts $2$-transitively on $\B_1'$ and $\B_2'$, we get that $\Gamma_{\B_1'}\cong \Gamma_{\B_2'}\cong K_p$. Since $D_j'\subset D$, we know that for any $g\in \overline{G}$, $d(X[B_{g}, B_{\tau g}])=2$. Let $C_{\B_1'}$ be a Hamilton cycle of $\Gamma_{\B_1'}$, which is

$$
B_{g_1},B_{g_2},B_{g_3},\ldots,B_{g_{p-1}},B_{g_p},B_{g_1}.
$$
Let $C_{\B_2'}$ be a Hamilton cycle of $\Gamma_{\B_2'}$, which is

$$
B_{\tau g_1},B_{\tau g_2},B_{\tau g_3},\ldots,B_{\tau g_{p-1}},B_{\tau g_p},B_{\tau g_1}. 
$$
Then
$$B_{g_1},B_{g_2},B_{g_3},\ldots,B_{g_{p-1}},B_{g_p}, B_{\tau g_p},B_{\tau g_{p-1}},\ldots,B_{\tau g_{2}},B_{\tau g_1},B_{g_1},
$$
denoted by $C_{\B'}$, is a Hamilton cycle of $\Gamma_{\B'}$. Since $d(X[B_{g_p}, B_{\tau g_p}])=2$, we can lift $C_{\B'}$ to a path of length $2p+1$, which is
$$
Hg_1',Hg_2',Hg_3',\ldots,Hg_{p-1}',Hg_p',H\tau \ell^{j}g_p',H\tau \ell^{-j}g_{p-1}',\ldots ,H\tau\ell^{j}g_1',H\ell^{2j} g_1',
$$
where $H\tau^i g_j'\in B_{\tau^ig_j}$. Then we can obtain $q$ vertex-disjoint paths of length $2p$, which are
$$
H\ell ^{kj}g_1',H\ell^{-kj}g_2',\ldots,H\ell ^{-kj}g_{p-1}',H\ell^{kj}g_p',H\tau \ell^{(k+1)j} g_p',H\tau\ell^{-(k+1)j}g_{p-1}',\ldots, H\tau\ell^{(k+1)j} g_1',
$$
where $k\in \ZZ_q$. In particular, $H\ell ^{kj}g_1'\sim H\tau\ell^{(k+1)j} g_1'$. After these $q$ paths are connected end to end, we obtain a Hamilton cycle in the graph.

Secondly, suppose $d\ge 3$. Without loss of generality, we may assume that $H<M<\overline{G}$, where $M$ is the point stabilizer of a 1-dimensional subspace of \(\mathbb{F}_{r}^d\). Set $$\ell=\begin{pmatrix}
\theta&\mathbf{0}^T&0\\
\mathbf{0}&E_{(d-2)\times {(d-2)}}&\mathbf{0}\\
0&\mathbf{0}^T&\theta^{-1}
\end{pmatrix}\in M\backslash H,$$
where $\mathbf{0}$ denotes the $(d-2)$-dimensional zero vector. Then $\ell^q\in H$. Set $g\in \overline{G}\backslash M$. By Corollary~\ref{co2.16}, we give all forms of the double cosets of $G$ with respect to $H$:

\begin{enumerate}[(i)]
	\itemsep=0pt

\item $H\tau^i\ell gH$, where $i\in \ZZ_2$, $2$ suborbits of length $q(p-1)$, both of them are self-paired.

\item $H\tau^i \ell^jH$, where $i\in\ZZ_2$, $j\in\ZZ_q$ and $(i,j)\neq (0,0)$, $2q-1$ suborbits of length $1$, all of them are not self-paired except for $H\tau H$.

\end{enumerate}

Since $\lg H,\ell,\tau \rg =N\times M\neq G$, we know that if $\Gamma=X(G,H,D)$ is a $G$-minimal graph, then $d(\Gamma)\ge q(p-1)\ge \frac{2}{3}pq$. By Proposition~\ref{pro2.1}, $\Gamma$ contains a Hamilton cycle. \qqed

\subsubsection{$\overline{G}$ acts non-quasiprimitively on $\B$}

\begin{theorem}
If $\overline{G}$ acts non-quasiprimitively on $\B$, then $\Gamma$ has a Hamilton cycle.
\end{theorem}

\demo In this case, let $\overline{T}$ be a maximal intransitive normal subgroup of $\overline{G}$. Then $\overline{T}$ induces blocks of length $r$, where $r\in \{p,q\}$. Let $\overline{P_0}$ be a Sylow-$r$ subgroup of $\overline{T}$, $\overline{T_0}=N_{\overline{G}}(\overline{P_0})$, and let $s=\frac{pq}{r}$. Then $\overline{P_0}\cong \ZZ_r^n$, where $n\le s$. Using the same argument as Lemma~\ref{le3.5}, we know that $\overline{T_0}$ acts transitively on $\B$ and $\overline{P_0}\triangleleft ~\overline{T_0}$. There exists $g_0\in \overline{T_0}$, where $o(g_0)=s^m$, such that $\overline{G_0}=\overline{P_0}\rtimes \lg g_0\rg$ acts transitively on $\B$. Thus, $G$ contains a subgroup $G_0=\ZZ_2.\overline{G_0}$, which is transitive on $V(\Gamma)$, where $|G_0|=2r^ns^m$. In particular, $\ZZ_2.\overline{P_0}\cong \ZZ_2\times \overline{P_0}$ and $\ZZ_2.\overline{P}_0\triangleleft G_0$. Let $P$ be the Sylow-$r$ subgroup of $\ZZ_2.\overline{P}_0$. Then $P \char N\times P\triangleleft G_0$ and so $P\triangleleft G_0$. Thus, $P$ induces blocks of length $r$. By Theorem~\ref{th:3.7}, $\Gamma$ contains a Hamilton cycle.\qqed

		\end{document}